\documentclass[english,11pt]{article}

\usepackage[english]{babel}
\usepackage{textcomp}
\usepackage{eurosym}
\usepackage{amsthm}
\usepackage[intlimits]{amsmath}
\usepackage{amssymb}
\usepackage{amsfonts}
\usepackage{bm}
\usepackage{parskip}
\usepackage{paralist}
\usepackage{graphics}
\usepackage{graphicx}

\usepackage{subfig}

\theoremstyle{plain}

\newtheorem{thms}{Theorem}[section]
\newtheorem{cors}[thms]{Corollary}

\newtheorem{props}[thms]{Proposition}

\theoremstyle{definition}
\newtheorem{defns}[thms]{Definition}
\newtheorem{rems}[thms]{Remark}
\newtheorem{exmps}[thms]{Example}

\numberwithin{equation}{section}

\renewcommand {\epsilon}{\varepsilon}

\newcommand{\EE}{\mathbb{E}}

\newcommand{\NN}{\mathbb{N}}

\newcommand{\PP}{\mathbb{P}}
\newcommand{\QQ}{\mathbb{Q}}
\newcommand{\RR}{\mathbb{R}}

\newcommand{\bB}{\mathcal{B}}
\newcommand{\cC}{\mathcal{C}}

\newcommand{\fF}{\mathcal{F}}

\newcommand{\lL}{\mathcal{L}}

\newcommand{\tT}{\mathcal{T}}

\newcommand{\al}{\alpha}

\newcommand{\e}{\varepsilon}
\newcommand{\la}{\lambda}

\newcommand{\si}{\sigma}

\newcommand{\om}{\omega}
\newcommand{\Om}{\Omega}

\newcommand{\ra}{\rightarrow}

\newcommand{\ti}{\tilde}

\newcommand{\ind}{\mathbf{1}}

\newcommand{\lqq}{\leqslant}
\newcommand{\gqq}{\geqslant}
\renewcommand{\le}{\leqslant}

\DeclareMathSymbol{\ophi}{\mathalpha}{letters}{"1E}

\renewcommand{\phi}{\varphi}

\newcommand{\be}{\begin{equation}}
\newcommand{\ee}{\end{equation}}
\newcommand{\ben}{\begin{equation*}}
\newcommand{\een}{\end{equation*}}

\newcommand{\ba}{\begin{equation}\begin{aligned}}
\newcommand{\ea}{\end{aligned}\end{equation}}

\DeclareMathOperator{\sign}{sign}

\newfont{\cyrfnt}{wncyr10}
\def\J3{\cyrfnt{\rm \u{\cyrfnt I}}}
\def\j3{\cyrfnt{\rm \u{\cyrfnt i}}}

\allowdisplaybreaks[4]

\usepackage[]{color}
\definecolor{DarkGreen}{rgb}{0.1,0.7,0.3}   

\begin{document}
\thispagestyle{empty}

\title{
Transportation distances and noise sensitivity \\
of multiplicative L\'evy SDE with applications
}

\author{Jan Gairing  \footnote{Institut f\"ur Mathematik, Humboldt-Universit\"{a}t zu Berlin, Berlin, Germany;  gairing@math.hu-berlin.de}
\hspace{1.5cm}
Michael H\"ogele\footnote{Departamento de Matem\'aticas, Universidad de los Andes, Bogot\'a, Colombia; ma.hoegele@uniandes.edu.co}
\hspace{1.5cm}
Tetiana Kosenkova \footnote{Institut f\"ur Mathematik, Universit\"at Potsdam, Potsdam, Germany; kosenkova@mail.math.uni-potsdam.de}
}

\maketitle

\begin{abstract}
This article assesses the distance between the laws of 
stochastic differential equations with multiplicative L\'evy noise 
on path space in terms of their characteristics. 
The notion of transportation distance 
on the set of L\'evy kernels 
introduced by Kosenkova and Kulik 
yields a natural and statistically tractable upper bound 
on the noise sensitivity. 
This extends recent results for the additive case 
in terms of coupling distances 
to the multiplicative case. 
The strength of this notion is shown in a statistical 
implementation for simulations and 
the example of a benchmark time series in paleoclimate. 
\end{abstract}


\section{Introduction}

Many dynamical phenomena are subject to random forcing, 
often described by stochastic differential equations of the following type
\begin{equation}\label{eq1}
dX(t) = -\nabla U(X(t)) dt + d \xi(t), \qquad X_0 = x_0,
\end{equation}
where $-\nabla U$ encodes the deterministic dynamics given by
a potential gradient and $\xi$ a noise signal.
In general, for instance when $\xi$ exhibits discontinuities, 
it is not straight-forward to describe the law of the solution $X$ on path space 
in terms of the parameters which determine the distribution of $\xi$. 
Our approach allows to quantify the distance of such laws in terms of accessible quantities, 
both analytically and statistically. 

The case where $\xi$ is given as a discontinuous L\'evy process 
was studied in a previous publication~\cite{GHKK14}.
The authors introduced the notion of a coupling distance
between L\'evy measures in order to quantify
the Wasserstein distance on path space. 
The coupling distances have been found to be sufficiently strong
(in a topological sense) to quantify the convergence in functional limit theorems,
yet being weak enough in order to be numerically and statistically tractable.
See for instance the calibration problem of a climate time series in~\cite{GHKK15}.
In many situations, however, the noise process $\xi$ exhibits state dependence, for instance multiplicative noise.  
This generalization lifts L\'evy diffusions to L\'evy-type diffusions. 
To treat this class of noise processes the authors introduced in~\cite{kos-kul}
the notion of transportation distance extending the coupling distances
with the help of a common reference L\'evy measure. 
The present article establishes analogous bounds on the distance between
the laws of L\'evy-type diffusions on path space in terms of transportation distances. 

\begin{figure}[p]
\begin{center}
\subfloat[ ]{\includegraphics[height=17em,page=1]{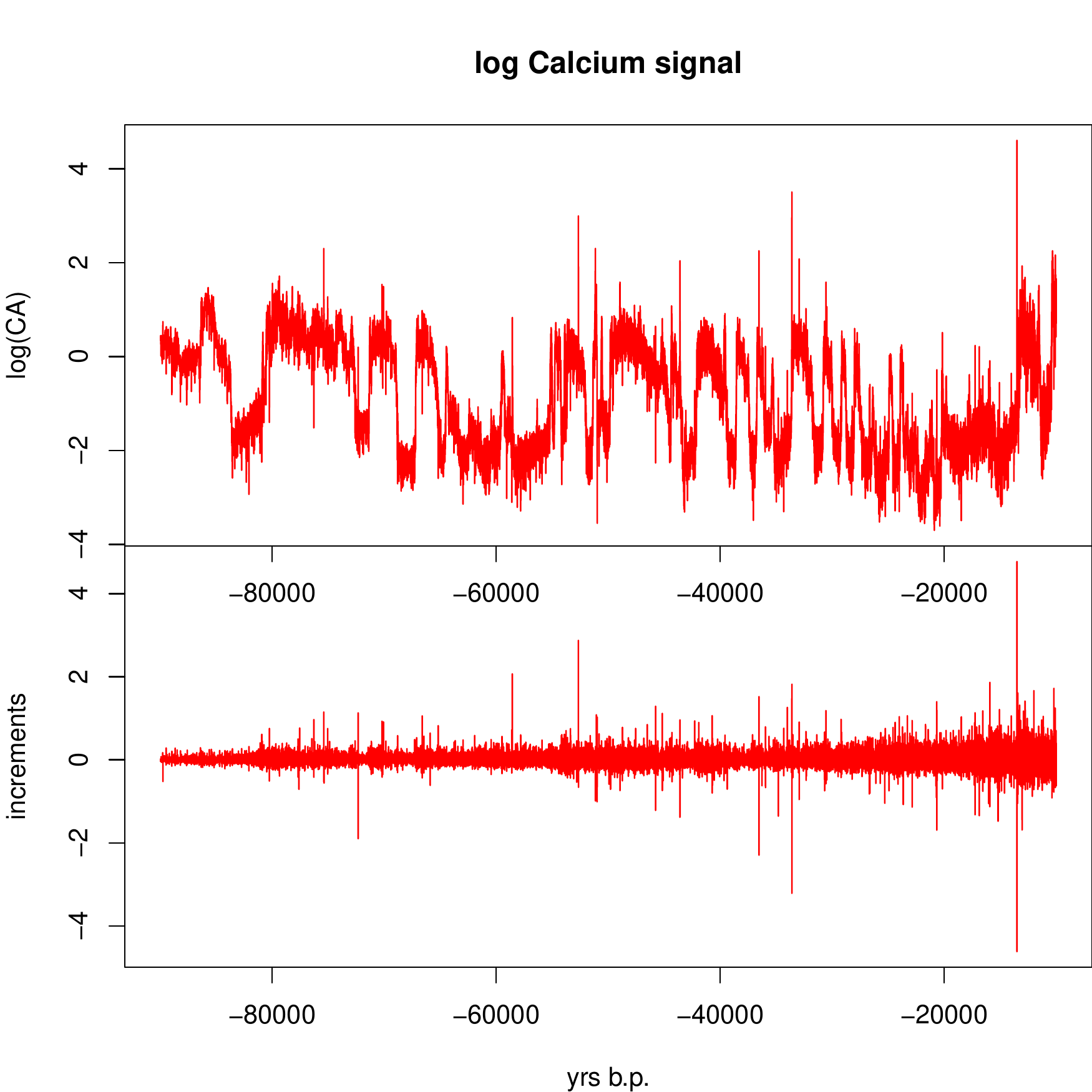}}
\end{center}
\caption{Timeseries}
\label{fig:TS}
\end{figure}

We stress that our procedure is suitable for a large variety of phenomena modelled with jump diffusions, 
such as in finance, e.g. \cite{Ra2003}, or neurosciences \cite{BerLan, DoThieu, TucRodWan}. 
The particular application we have in mind in this article is the refinement of 
the analysis of the noise structure behind the paleoclimate temperature evolution studied in \cite{GHKK15}. 
The climate data apparently fluctuate around two distinct metastable states 
with rapid transitions (see Fig. \ref{fig:TS}). 
Such phenomena are observed in 
stochastic energy balance models \cite{Arn, BenziEtal, BerGen, Ha76, Im01, IM02}. 
There is a list of publications associating this time series to 
an underlying jump diffusion, see for instance \cite{Dit99a,GHKK15,GI14,HIP09}.
Using various techniques these articles aim to determine 
the (polynomial) jump behavior of \eqref{eq1} 
for different classes of heavy-tailed L\'evy processes $\xi$.
The models so far require 
the spacial homogeneity of the noise characteristics. 
The present article lifts this restriction.
We may now investigate the statical behavior in the different 
spatial regimes, prescribed by the metastable states,
and solve the corresponding model selection problem on the generic class of heavy-tailed jump diffusions. 
Our results of the implementation of this program applied to the mentioned climate times series 
are consistent with the findings in \cite{GHKK15}.

\section{Preliminaries}

\paragraph{Transportation functions and transportation distance: }
Consider a filtered probability space $(\Om, \fF, \PP, (\fF_t)_{t\gqq 0})$ satisfying the 
usual conditions in the sense of Protter \cite{Pr04} carrying a scalar Brownian motion $(B_t)_{t\in [0, T]}$ and
an independent Cauchy Poisson random measure $\nu_0$ on $[0, T] \times \RR$ with intensity measure $dt \otimes \Pi_0$ 
given by
$$\Pi_{0}(dv)=\ind_{\{v\in \RR\setminus \{0\}\}}\frac{dv}{v^{2}}.$$
We define a L\'evy measure to be a $\si$-finite Borel measure on $\RR\setminus\{ 0 \}$ satisfying 
\begin{equation}
\int_{\RR}(|v|^2\wedge 1)\Pi(dv)<\infty.
\label{eq:levymeasure}
\end{equation}
In contrast to the standard definition we do not exclude point-mass in $0$ which may be taken to be infinity.
Nevertheless we will identify all such measures that coincide on the Borel $\sigma$-algebra $\bB(\RR\setminus\{0\})$ 
and the standard L\'evy measures (without mass in $0$) 
are the canonical representatives of those equivalence classes (for details see \cite{kos-kul}).
The key of our analysis is the observation that all standard L\'evy measures with infinite mass admit a
representation as a measure transform of a common reference L\'evy measure.
This reference measure will be chosen to be the standard symmetric Cauchy measure.
In order to treat finite L\'evy measures we may artificially assign an infinite point-mass to $0$.
The precise statement is given as follows.

\begin{props}\label{prop:transform}
For any standard L\'evy measure $\Pi$ (that is $\Pi(\{0\}) = 0$) 
there exists a unique L\'evy measure $\ti \Pi$ in the sense of (\ref{eq:levymeasure}) defined as 
$\ti \Pi = \Pi$ if $\Pi(\RR) = \infty$ and $\ti \Pi = \Pi + \infty\cdot\delta_0$ else,  
and a unique \emph{transportation function} $c:\RR\ra\RR$ satisfying 
    \begin{enumerate}
     \item $c$ is non-decreasing, 
    \item  $c((-\infty,0]) \subset (-\infty,0]$ and $c([0, \infty)) \subset [0, \infty)$,
    \item $c$ is left continuous on $(-\infty, 0]$ and right continuous on $[0,\infty)$,
    \end{enumerate}
such that $\ti \Pi=\Pi_{0}\circ c^{-1}$.
\end{props}
\noindent A proof is found in \cite{kos-kul}. 
For the sake of readability and due to uniqueness 
we no longer distinguish between $\ti \Pi$ and $\Pi$.

\begin{exmps}\label{ex:powerlaw}
\begin{enumerate}
 \item For the one-sided Cauchy measure $\Pi(dv)=\mathbf{1}_{\{v>0\}}\frac{dv}{v^{2}}$ 
 we have the transportion function $c(v)=\mathbf{1}_{\{v>0\}}v.$
Note that $c$ maps the infinite mass of $\Pi_0$ on $(-\infty,0]$ to the point $0$.
The image measure $\tilde{\Pi}(dv)=\mathbf{1}_{\{v>0\}}\frac{dv}{v^{2}}+\infty\cdot\delta_{0}$ coincides with $\Pi$ outside 0.
 \item
 The next example will be exploited extensively in the applications of Section \ref{sec:application}.
 For $\alpha,\e,\lambda>0$ consider the Pareto-type power tail $\Pi_{\alpha,\varepsilon,\lambda}(dv) = \ind_{\{v>\e\}} \frac{\lambda dv }{v^{\al+1}}$.
 It is easy to verify that in this case the transportation function is of the form 
 \[
 c(v) = \ind_{\{v>\tfrac{\alpha \varepsilon^\alpha}{\lambda} \}} \,\left(\frac{\lambda v}{\alpha} \right)^{1/\al}  \ .
 \]
 Note that the positive span of the family $\left\{  \Pi_{\alpha,\varepsilon,\lambda}: \alpha,\varepsilon,\lambda\in \QQ_+ \right\}$ 
 is dense in the class of heavy tailed L\'evy measures on the positive half-line.
\end{enumerate}
\end{exmps}

\begin{exmps}\label{ex:emprirical}
Let $(\xi_i)_{i\in \NN}$ be a sequence of real valued random variables with values in $[\e, \infty)$ for some $\e>0$.
The respective \emph{empirical measure}  is
given by $\Pi_n([a,b]) = \frac{\#\{\xi_i \in (a,b], i\le n\}}{n}$ on intervals $(a,b]$. 
We denote the $n$-th \emph{order statistic} by $\e< \xi_{1:n}\le \xi_{2:n}\le\cdots\le\xi_{n:n}$. 
In each interval $(\xi_{i:n}, \xi_{i+1:n}]$ we have exactly one data point with individual mass $1/n$ such that 
\[
\Pi_n((\xi_{i:n}, \xi_{i+1: n}]) = \frac{1}{n} = \int_{c^{-1}(\xi_{i:n})}^{c^{-1}(\xi_{i+1: n})} \frac{dv}{v^2} 
= \frac{1}{c^{-1}(\xi_{i:n})}- \frac{1}{c^{-1}(\xi_{i+1:n})},
\]
which yields
\begin{align*}
c^{-1}(v) =
\begin{cases}
0 & v\in [0, \e) \\
1 & v\in [\e, \xi_{1:n})\\
\frac{n}{(n-i+1)} & v\in [\xi_{i-1:n}, \xi_{i:n})\qquad \mbox{ for all }i\in \{2, n\}\\
\infty & v \gqq \xi_{n:n}\ .
\end{cases}
\end{align*}
Hence the inverse is given by 
\begin{align*}
c(v) =
\begin{cases}
0 & v \in  [ 0,1 ) \\
\xi_{i:n} & v\in [\frac{n}{(n-i+1)},\frac{n}{(n-i)} )\qquad \mbox{ for all }i\in \{1, n-1\}\\
\xi_{n:n} & v \gqq n.
\end{cases}
\end{align*}
\end{exmps}

\noindent We take advantage of the measure transform of Proposition \ref{prop:transform} to compare two given L\'evy measures in 
terms of a truncated $L^p(\Pi_0)$ distance of the respective transportation functions.

\begin{rems}\label{rem:allLevy} We stress that Proposition \ref{prop:transform} guarantees that 
the previously chosen probability space $(\Om, \fF, \PP, (\fF_t)_{t\gqq 0})$ 
is rich enough to carry any Poisson random measures $\nu$ (in distribution) 
with respect to the intensity measure $dt \otimes \Pi$, where $\Pi$ is a (standard) 
L\'evy measure via 
\[
\nu((s, t] \times A) = \nu_0((s, t] \times c^{-1}(A)), \qquad A\in \bB(\RR\setminus \{0\}), \quad 0\lqq s \lqq t.  
\]
\end{rems}

\begin{defns}\label{tm}
Define the class
\[
\lL^p_\rho := \{\Pi \mbox{ (standard) L\'evy measure}~|~\int_\RR \rho^p(0, v) \Pi(dv)<\infty\},
\]
where $p\gqq 1$ and $\rho(x,y) := |x-y|\wedge 1$ is the truncated Euclidean distance.
For L\'evy measures $\Pi_{i}\in \lL^p_\rho$, $i=1,2$ 
with respective representations $\Pi_i = \Pi_0 \circ c_i^{-1}$
according to Proposition \ref{prop:transform} 
we define the \textit{transportation distance} of order $p$ by
\begin{align*}
\mathcal T_p(\Pi_{1},\Pi_{2}):=\bigg( \int_{\RR}\rho^{p}\bigl(c_1(u),c_2(u)\bigr)\,\Pi_{0}(dv) \bigg)^{1/p}.
\end{align*}
\end{defns}
\noindent It has been established in \cite{kos-kul}
that the transportation distance metrizes 
the positive cone~$\lL^p_\rho$.

\begin{rems}
In the case of proper L\'evy-type mesures
with state dependence $\Pi(x, dv)$
Proposition~\ref{prop:transform} yields a family of
transportation functions $(c(x, \cdot))_{x\in \RR}$.
\end{rems}

\paragraph{Stochastic differential equations: }

For any Poisson random measure $\nu$ on $(\Omega,\fF,\PP, (\fF_t)_{t\gqq 0})$ (cf. Remark \ref{rem:allLevy}) 
independent of $W$ with intensity measure $dt \otimes \Pi$, $\Pi$ a L\'evy measure, 
and $a,b\in\RR$ there is a L\'evy process $Z$ given in terms of its L\'evy-It\^o decomposition  
\begin{equation}
Z(t) = at + b W(t) + \iint_{(0,t]\times [-1,1]} v\tilde\nu(dv,dt) + \iint_{(0,t]\times \RR\setminus [-1,1]} v\nu(dv,dt).
\label{eq:levyito}
\end{equation}
As usual $\ti \nu = \nu - dt \otimes \Pi$ denotes 
the compensated Poisson random measure of $\nu$ on $\bB(\RR \setminus\{0\})$. 
Recall that the law of $Z$ is characterized by its \emph{L\'evy triplet} $(a,b,\Pi)$. 
For further details we refer to \cite{Sato-99}. 
Consider now the following formal stochastic differential equation
\begin{align}
\begin{aligned}
dX(t)
&= a(X(t)) dt + b(X(t)) dW(t) \\
&\quad + \int_{|v|\lqq 1} \,v\,  \ti \nu(X(t-), dt, dv) + \int_{|u|> 1} \,v\, \nu(X(t-), dt, dv). \label{eq_ij}
\end{aligned}
\end{align}
Here $\nu$ (resp. $\ti \nu$) is interpreted as a space dependent (compensated) Poisson random measure explained below. 
A solution $X$ of such an equation is understood as the solution to 
the corresponding martingale problem for the following integro-differential operator $A$ 
acting on $\phi\in \cC^2_c(\RR)$
\begin{align*}
A[\phi](x) 
&= a(x) \phi'(x) + b^2(x) \phi''(x) + \int_{\RR\setminus \{0\}} \big(\phi(x+v) - \phi(x) - \phi'(x) v \ind_{\{|v|\lqq 1\}}\,\big) \Pi(x, dv),
\end{align*}
where $a, b: \RR\ra \RR$ and $x\mapsto \Pi(x, \cdot)$ is a L\'evy kernel, 
which associates to each $x\in \RR$ the L\'evy measure $\Pi(x, \cdot)$. 
We may rewrite $A$ in terms of the Lipschitz continuous cutoff function $\tau(u) := \sign(u) (|u|\wedge 1)$ 
and $\bar a(x) = \Pi(x, \{u\in\RR~|~|u|>1\})$ as 
\begin{equation}\label{eq:martingale problem}
A[\phi](x) 
= (a(x) +\bar a(x)) \phi'(x) + b^2(x) \phi''(x) + \int_{\RR\setminus \{0\}} \big(\phi(x+v) - \phi(x) - \phi'(x) \tau(v)\big) \Pi(x, dv).
\end{equation}
Proposition \ref{prop:transform} allows us to represent $\Pi(x, \cdot)$
as $\Pi_0 \circ c^{-1}(x, \cdot)$ in terms of a family of tansport functions $(c(x, \cdot))_{x\in \RR}$
with respect to the Cauchy reference measure $\Pi_0$. For L\'evy type processes we refer to 
\cite{GS82} and \cite{BSW13}.

In \cite{kos-kul} it is shown that under the following Lipschitz and boundedness conditions on the space-dependent
coefficients $(a, b, \Pi)$ for any $x, y\in \RR$
\begin{align}
&|a(x)-a(y)| \lqq L_{a} \,|x-y|,\label{Lip0a}\\[2mm]
&(b(x)-b(y))^{2} \lqq L_{b} \,\rho^{2}(x,y),\label{Lip0b}\\[2mm]
&|\bar a(x) - \bar a(y)| \lqq \,L_{\bar a}\rho(x,y),\label{Lip0e}\\[2mm]
&\mathcal T^{2}_2(\Pi(x,\cdot),\Pi(y,\cdot)) \lqq L_{\Pi} \,\rho^{2}(x,y),\label{Lip0c}
\end{align}	
and initial value $x_0 \in \RR$ there exists a unique strong solution to the 
martingale problem associated to (\ref{eq:martingale problem}). 
On the probability space $(\Om, \fF, \PP, (\fF_t)_{t\gqq 0})$ 
this is given as a strong solution of the following SDE
\begin{align}
\label{eq:sde2}
dX(t)
&= \big(a(X(t)) + \bar a(X(t))\big)dt + b(X(t)) dW(t) \nonumber\\
&\quad + \int_\RR c(X(t-), v) \left[  \nu_0(dt, dv)-\frac{\tau(c(X(t-), v))}{c(X(t-), v)}\Pi_{0}(dv)dt\right], \\
X(0) &= x_0.
\end{align}
Under these assumptions Remark \ref{rem:allLevy} 
provides a strong solution of (\ref{eq_ij}). 

The main purpose of this construction is 
to compare the laws of two diffusions $X_i$, $i=1,2$ 
on the same probability space $(\Om, \fF, \PP)$
given as strong solutions (\ref{eq:sde2})
with identical Brownian motion $W$ and Cauchy Poisson random measure $\nu_0$,
whose respective coefficients $(a_i, b_i, \Pi_i)$, $i=1, 2$ 
satisfy (\ref{Lip0a}), (\ref{Lip0b}), (\ref{Lip0e}) and (\ref{Lip0c}).
This is carried out in Theorem \ref{thm:1}. 
In a simplified setting more suitable for the applications 
we have in mind this is done in Theorem \ref{thm_estim_for_1_rho}.

\begin{rems}\label{rem: T1}
In the case of pure jump diffusions with finite intensity, 
that is L\'evy-type diffusions with triplet of characteristics $(a, 0, \Pi)$, 
where $\Pi$ is a finite measure, it is obvious that 
we only need the Lipschitz continuity of $a$ (\ref{Lip0a}). 
In this case we consider equation 
\begin{align}
\label{sde3}
dX(t)
&=  a(X(t))dt + \int_\RR c(X(t-), v) \nu_0(dt, dv), \\
X(0) &= x_0.
\end{align}
\end{rems}

\section{Main results}

In this section we will first compare the laws of two diffusions in the sense of (\ref{eq:sde2})
in terms of $\tT_2$. In a second part this is carried out in terms of $\tT_1$
motivated by applications elaborated in Section~\ref{sec:application}.

\subsection{Noise sensitivity estimates in terms of $\tT_2$}

We are interested in the sensititivity of the laws of solutions of (\ref{eq:sde2}) 
with respect to their parameters. The following theorem provides a quantitative estimate 
in terms of $\tT_2$. 

\begin{thms}\label{est_forT_2_rho}
Let $X_i$, $i= 1,2$ be strong solutions of (\ref{eq:sde2}) on $(\Om, \fF, \PP)$ with respective initial conditions $x_i\in \RR$
and triplets of characteristics $(a_i, b_i, \Pi_i)$ satisfying the Lipschitz conditions (\ref{Lip0a}-\ref{Lip0c})
for common constants $L_{a}, L_{b}, L_{\bar a}, L_{\Pi}>0$ and $\|a_1 - a_2\|_\infty < \infty$. 
Then for any $T>0$ there is constant $K>0$ such that for $G(x) = \max\{\sqrt{x}, x\}, x\gqq 0$ we have
\ba\label{thm_estim}
\EE\sup_{t \in [0, T]} \rho^2(X_1(t),X_2(t))
&\lqq K G\big(\Delta\big),
\ea
where 
\begin{align*}
\Delta &= \rho(x_1, x_2) + \|a_1-a_2\|_\infty^2 +\|\bar a_1-\bar a_2\|_\infty^2+ \|b_1 - b_2\|_{\infty}^{2} \\[2mm]
&\qquad + \sup_{x\in\RR}\mathcal T_2(\Pi_1(x, \cdot), \Pi_2(x, \cdot))
+\sup_{x\in\RR}\mathcal T_2^{2}\big(\Pi_1(x, \cdot), \Pi_2(x, \cdot)\big). 
\end{align*}
\end{thms}

The statement of the theorem also allows to compare laws of solutions 
of (\ref{eq_ij}) which are not necessarily defined on the same probability space. 
The natural choice of a metric between such laws is the (truncated) 
Wasserstein distance on path space, for details see \cite{GHKK14} 
and for general reference~\cite{RR98}. 

\begin{cors}
Under the assumptions of Theorem \ref{est_forT_2_rho} 
the strong solutions $X_i$, $i=1,2$ to (\ref{eq:sde2}) 
form a specific coupling of solutions to (\ref{eq_ij}). 
Hence the bound (\ref{thm_estim}) gives an upper bound of the 
Wasserstein distance of order $2$ on the path space of c\`adl\`ag paths endowed 
with the (truncated) supremum norm.  
The Wasserstein distance is precisely defined as 
as the infimum of the left-hand side of (\ref{thm_estim}) 
where the pair $(X_1, X_2)$ ranges over all couplings of $X_1$ with $X_2$. 
\end{cors}

The following example shows that the Lipschitz conditions (\ref{Lip0a}-\ref{Lip0c}) 
are not very restrictive in the class of L\'evy-type kernels with exponential moments.

\begin{exmps}
An important example is given by the Gamma-type process, that is a L\'evy-type process
with triplet of characteristics $(0, 0, \Pi)$ and initial value $x_0 = 0$.
The L\'evy kernel is given as $\Pi(x,dv)=\ind_{\{v> 0\}}\lambda(x)(e^{-\gamma(x)v} /v) dv$ with bounded and Lipschitz continuous coefficients
$\la, \gamma: [0, \infty) \ra [0, \infty)$ satisfying $\sup_{x\in\RR}\lambda(x)>0$ and $\sup_{x\in\RR}\gamma(x)>0.$
Indeed, Proposition 3 in \cite{kos-kul} states the following.
Let $\Pi_j$, $j=1,2$ be two Gamma L\'evy measures
$\Pi_j(dv)=\ind_{\{v>0\}}(\gamma_j e^{-\lambda_j v}/v)dv$, $j=1,2$ with constants $\ \gamma_j,\lambda_j>0$.
\begin{enumerate}
\item\label{gamma} For two such Gamma measures $\Pi_j$, $j=1,2$
with the same parameter $\lambda$ and different parameters $0<\gamma_{1} < \gamma_{2}$
there exists a constant $D = D(\la) >0$ such that the following bound holds true
\ba\label{T_bound}
\mathcal T_2(\Pi_{1},\Pi_{2}) \lqq D\left(\gamma_{2}-\gamma_{1}\right).
\ea
\item\label{lambda}
For two such Gamma measures $\Pi_j$, $j=1,2$
with the same parameter $\gamma$ and different parameters $0<\lambda_{1} < \lambda_{2}$
there exists a constant $\tilde{D} = \ti D(\gamma, \la_1, \la_2) >0$ locally bounded in $\la_1$ 
around $\la_2$ such that 
\ba\label{T_bound_2}
\mathcal T_2(\Pi_{1},\Pi_{2}) \lqq \tilde{D}\left(\lambda_{2}-\lambda_{1}\right).
\ea
\end{enumerate}
These results immediately yield that condition (\ref{Lip0c}) is satisfied in either case.
 \end{exmps}

\subsection{Noise sensitivity estimates in terms of $\tT_1$}

In applications, however, the involved L\'evy measures are heavy-tailed and do not satisfy (\ref{Lip0c}).
In \cite{GHKK15} for instance, it was obtained that the measures that are likely to describe the noise
in the paleoclimatic data are of the following shape
 \ba\label{pol_meas}
\Pi(dv)&= \Big(\mathbf{1}_{\{v<1\}}\frac{\lambda^-}{|v|^{\alpha^-+1}} + \mathbf{1}_{\{v>1\}} \frac{\lambda^+ }{v^{\alpha^++1}}\Big) dv,\quad \alpha\gqq 2.
\ea
\begin{rems}
 Such L\'evy measures do not satisfy (\ref{Lip0c}) due to the following reasoning. 
Let $\la^+=1$, $\la^- = 0$ and $\al_2 > \al_1 > 1$ with $|\al_2-\al_1|$ small. 
Then transportation distance 
satisfies inequalities of the form 
\begin{align*}
\tT_2^2(\Pi_1, \Pi_2) 
&= \int_0^\infty \Big(\ind_{\{v>\al_1^{1/\al_1}\}} ~\Big(\frac{v}{\al_1}\Big)^{1/\al_1}
- \ind_{\{v> \al_2^{1/\al_2}\}} ~\Big(\frac{v}{\al_2}\Big)^{1/\al_2}\Big)^2 \wedge 1)  \Pi_0(dv)  \\
&\gqq \int_{\al_1^{1/\al_1}}^{\al_2^{1/\al_2}} \Big(\frac{v}{\al_1}\Big)^{1/\al_1}\frac{1}{v^2} dv 
\gqq C_{\al_1} |\al_2^{1/\al_2}- \al_1^{1/\al_1}| \gqq \inf_{|x|\lqq |\al_2-\al_1|} g'(x) ~C_{\al_1} |\al_1- \al_2|,
\end{align*}
where $g(x) = x^{1/x}$. Since $g'(\al_1)$ is locally bounded from below in a small neighborhood of $\al_1$,  
this only implies H\"older continuity with respect to $\tT_2$. 
\end{rems}

\noindent With this result in mind we switch to the setting of a pure jump diffusions $X_j$
with the respective triplet of characteristics $(a_j, 0, \Pi_j)$ for finite $\Pi_j$ 
and initial conditions $x_j$, satisfying (\ref{Lip0a}).
For two solutions of such differential equations we can state the following result.

\begin{thms}\label{thm:1}
Let $X_i$, $i =1,2$ be strong solutions of (\ref{eq:sde2})  on $(\Om, \fF, \PP)$ 
with respective initial conditions $x_j \in \RR$ and triplet of characteristics $(a_i, 0, \Pi_i)$ 
satisfying the Lipschitz conditions
(\ref{Lip0a}) and 
\begin{equation}\label{Lip0d}
\tT_1(\Pi(x, \cdot), \Pi(y, \cdot)) \lqq L_\Pi \rho(x,y) \qquad \mbox{ for all }x,y\in \RR.
\end{equation}
Then for any $T>0$ there exists a constant $K>0$ such that for $G(x) = \max\{\sqrt{x}, x\}, x\gqq 0$ we have
\ba\label{thm_estim_for_1_rho}
&\EE\sup_{t \in [0, T]}\rho(X_1(t), X_2(t))\lqq K G\Big(\rho(x_1, x_2) + \sup_{x\in\RR}\mathcal T_1\big(\Pi_{1}(x,\cdot),\Pi_{2}(x,\cdot)\big)\Big).\\
\ea
\end{thms}
\noindent The following proposition establishes that for constants $\la^+\gqq 0$ and $\la^-\gqq 0$ and
a bounded and globally Lipschitz continuous function $\al: \RR \ra [2, \infty)$ the L\'evy-type measure $\Pi(x, dv)$
\ba\label{pol_meas_type}
\Pi(x, dv)&= \Big(\mathbf{1}_{\{v<1\}} ~\frac{\lambda^-}{|v|^{\alpha(x)+1}}+ \mathbf{1}_{\{v>1\}} ~\frac{\lambda^+ }{v^{\alpha(x)+1}}\Big) dv
\ea
satisfies the $\tT_1$-Lipschitz condition (\ref{Lip0d}).

\begin{props}\label{prop:1}
For L\'evy measures $\Pi_{i}$, $i =1,2$ defined in (\ref{pol_meas})
with $\la^+ = \la_1^+ = \la_2^+$ and $\la^- = \la_1^- = \la_2^-$
and
$2 \lqq \al_1^+, \al_2^+, \al_1^-, \al_2^-$
there exists a constant $D = D(\la^+, \la^-) >0$ such that
\ba
\mathcal T_1(\Pi_{1},\Pi_{2})\lqq D\Big(|\alpha_{1}^+-\alpha_{2}^+|+|\alpha_{1}^--\alpha_{2}^-|\Big).
\ea
\end{props}
\noindent All proofs are found in Section \ref{sec: proofs}.

\section{Applications}
\label{sec:application}

In this section we demonstrate the benefit of the transportation distance in an empirical setting.
Assume we are given a data set $\xi_1,\ldots \xi_n$ that we interpret as jumps of a L\'evy process with L\'evy measure $\Pi$.
In Example \ref{ex:emprirical} we have calculated the transportation kernel $c_n$ of the empirical L\'evy measure $\Pi_n$.
We are particularly interested in modelling with power law tails where the corresponding 
transportation kernels $c_\alpha$ are given in Example \ref{ex:powerlaw}. We are then in the 
position to evaluate the transportation distance between the empirical L\'evy measures $\Pi_n$ and 
such power laws. This is done in example \ref{ex:empiricaldisance}.
In the sequel we analyze the behavior of the transportation distance for simulated power law jumps.
Ultimately we utilize this device to propose a simple but state dependent jump diffusion model for a paleoclimatic bench mark time series.

\subsection{Analytical considerations}

\begin{exmps}
\label{ex:empiricaldisance}
Let us consider the transportation distance between the empirical measure $\Pi_n((a,b])$ 
of Example \ref{ex:emprirical} and a power law $\Pi_{\alpha,\e}(vu) =  \ind_{\{v>\e\}} \frac{\lambda dv }{v^{\al+1}} $
($\alpha\ne 1$) of Example \ref{ex:powerlaw}.
Recall that $c_n$ is piecewise constant. Observe that for $0<a<b$ and a constant $c>0$
\begin{align}
\int_a^b |c-(\tfrac{\lambda v}{\alpha} )^{1/\alpha}| v^{-2} dv
&= \int_{a\wedge \frac{\alpha}{\lambda} c^{\alpha}}^{b \wedge \frac{\alpha}{\lambda} c^{\alpha}} (c-(\tfrac{\lambda v}{\alpha} )^{1/\alpha} ) v^{-2} dv
+ \int_{a \vee \frac{\alpha}{\lambda} c^{\alpha}}^{b\vee\frac{\alpha}{\lambda} c^{\alpha}}  ((\tfrac{\lambda v}{\alpha} )^{1/\alpha}  -c) v^{-2} dv\nonumber\\
&=
q_c(b\vee\tfrac{\alpha}{\lambda} c^{\alpha}) - q_c(b\wedge\tfrac{\alpha}{\lambda} c^{\alpha}) +q_c(a\wedge\tfrac{\alpha}{\lambda} c^{\alpha}) - q_c(a\vee\tfrac{\alpha}{\lambda}c^{\alpha}),
\label{eq:L2}
\end{align}
where
\begin{equation*}
q_c (x) = (\tfrac{\lambda}{\alpha})^{1/\alpha}\frac{x^{1/\alpha-1}}{1/\alpha-1} +cx^{-1}\ .
\end{equation*}
For convenience we denote
\begin{equation*}
\ti \kappa_c^- := \tfrac{\alpha}{\lambda}(c-1)^{\alpha} ~\lqq~
\ti \kappa_c^0 := \tfrac{\alpha}{\lambda}(c)^{\alpha} ~\lqq~
\ti \kappa_c^+ := \tfrac{\alpha}{\lambda}(c+1)^{\alpha} \ .
\end{equation*}
Let us restrict these quantities to the interval $(a,b]$ by setting
\begin{equation*}
\kappa_c^\star:= \kappa_c^\star(a,b) := (a\vee \ti \kappa_c^\star )\wedge b \ , \quad \star\in \left\{ -,0,+  \right\} \ ,
\end{equation*}
which yields $a~\lqq~ \kappa_c^-  ~\lqq~
\kappa_c^0  ~\lqq~
\kappa_c^+ ~ \lqq~ b$.
Then formula (\ref{eq:L2}) turns into 
\begin{align*}
\label{eq:problem}
&\int_a^b (|c-(\tfrac{\lambda v}{\alpha} )^{1/\alpha}|\wedge 1) v^{-2} dv\\
&=\int_{a }^{\kappa_c^-}  v^{-2} dv
+ \int_{\kappa_c^-}^{\kappa_c^0} (c-(\tfrac{\lambda v}{\alpha} )^{1/\alpha}) v^{-2} dv
 + \int_{\kappa_c^0}^{\kappa_c^+} ((\tfrac{\lambda v}{\alpha})^{1/\alpha}-c) v^{-2} dv
 + \int_{\kappa_c^+}^b  v^{-2} dv \\
&=  q_c\bigl( \kappa_c^- \bigr)
 +	 q_c\bigl( \kappa_c^+ \bigr)
 - 2q_c\bigl( \kappa_c^0 \bigr) 
 +  \frac{1}{a} - \frac{1}{b} + \frac{1}{\kappa_c^+} - \frac{1}{\kappa_c^-}.
\end{align*}
Hence we may evaluate 
\begin{align*}
\mathcal T_1 (\Pi_n,\Pi_{\alpha, \e})
&= \int_0^{\infty} (|c_n(v) - c_\alpha(v)|\wedge 1) v^{-2} dv\\
&=  \int_\varepsilon^1 ((\tfrac{\lambda v}{\alpha} )^{1/\alpha}\wedge 1) v^{-2} du 
+ \sum_{i=1}^n \int_{\frac{n}{n-i+1}}^{\frac{n}{n-i}} (|\xi_{i:n} -(\tfrac{\lambda v}{\alpha} )^{1/\alpha}|\wedge 1) v^{-2} dv\\
&= \sum_{i=1}^n \int_{\frac{n}{n-i+1}}^{\frac{n}{n-i}} (|\xi_{i:n} -(\tfrac{\lambda v}{\alpha} )^{1/\alpha}|\wedge 1) v^{-2} dv\\
&= 1+\sum_{i=1}^n  q_{\xi_{i:n}}\bigl( \kappa^-(i) \bigr)
 +	 q_{\xi_{i:n}}\bigl( \kappa^+(i) \bigr)
 - 2q_{\xi_{i:n}}\bigl( \kappa^0(i) \bigr)
 + \frac{1}{\kappa^+(i)} - \frac{1}{\kappa^-(i)},
\end{align*}
where we have abbreviated
\begin{align*}
 \kappa^\star(i):= \kappa_{\xi_{i:n}}^\star\bigl(\tfrac{n}{n-i+1},\tfrac{n}{n-i}\bigr)  \ , \quad i=1,\ldots,n,\ \star\in \left\{ -,0,+  \right\} \ .\\
\end{align*}
\end{exmps}

\begin{rems}[Normalization]
\hfill
\begin{enumerate}
\item
Note that the empirical measure $\Pi_n$ is a probability measure and has total mass one.
It is therefore reasonable to normalize $\Pi_{\alpha,\e}$ choosing $\lambda=\frac{\varepsilon^{\alpha}}{\al}$.
\item In general the transportation distance is unbounded. In this example
the support of the L\'evy measures it bounded away from zero by $\varepsilon$. Hence the transportation distance 
is bounded from above by $\int_\varepsilon^\infty u^{-2}du = \varepsilon^{-1}$.
For the sake comparability we therefore normalize $\mathcal T_1(\Pi_n,\Pi_{\alpha,\e})$ and set
\begin{equation}
\tilde{\mathcal T}_1(\Pi_n,\Pi_{\alpha,\e}) := \varepsilon \cdot\mathcal  T_1(\Pi_n,\Pi_{\alpha,\e}) \ .
\label{eq:Ttilde}
\end{equation}
\end{enumerate}
\end{rems}

\subsection{Simulated data}\label{sec:simulation}

Let us assess the statistical behavior of the transportation distance between power laws and empirical counterparts in a simulation study.
To be precise we consider a family of L\'evy measures $(\Pi_{\alpha,\varepsilon})_{\alpha,\varepsilon}$ 
defined in Example \ref{ex:powerlaw}.2, where the intensity parameter $\lambda=\frac{\varepsilon^{\alpha}}{\al}$ 
is chosen such that it is a probability measure ($\Pi_{\alpha, \e}(\RR_+)=1$).
For $\alpha,\varepsilon$ fixed, we generate an i.i.d. sample $\xi_1,\ldots,\xi_n$ (for $n=100$) 
distributed according to $\Pi_{\alpha,\varepsilon}$, and evaluate the transportation distance $ \tilde{\mathcal T}_1$ 
between its empirical measure $\Pi_n$ and $\Pi_{\alpha,\varepsilon}$.
This experiment is repeated 100 times and the mean value is given in Table \ref{fig:mean100}, the corresponding standard deviation in Table \ref{fig:sd100}.
It is apparent that the magnitude of the renormalized distance $\tilde{\mathcal T}_1(\Pi_n,\Pi_{\alpha,\varepsilon})$ changes considerably over the parameter range.
However the standard deviation is always smaller than mean value by an order of magnitude.
Increasing $\alpha$ from $1$ to $10$ results in a decrease of the mean distance by more than a digit.
However increasing $\varepsilon$ tends to increase $\mathcal T_1$ grosso modo.

This behavior is explained by the Cauchy weighting with $v^{-2}$. Increasing $\alpha$ shifts mass away from zero to larger values which are damped stronger.

\begin{table}
\begin{center}
\begin{tabular}{rr|cccccc}
&$ \varepsilon$   & 0.5    & 0.6    & 0.7    & 0.8    & 0.9    & 1.0   \\[-3mm]
$\alpha$ &	&	&	&	&	&	&	\\
\hline
1  && 0.1901 & 0.2000 & 0.2191 & 0.2429 & 0.2766 & 0.3191\\
2  && 0.0877 & 0.0893 & 0.1039 & 0.1017 & 0.1181 & 0.1406\\
3  && 0.0522 & 0.0530 & 0.0579 & 0.0619 & 0.0690 & 0.0789\\
4  && 0.0364 & 0.0359 & 0.0373 & 0.0400 & 0.0479 & 0.0551\\
5  && 0.0286 & 0.0277 & 0.0272 & 0.0294 & 0.0335 & 0.0403\\
6  && 0.0232 & 0.0213 & 0.0216 & 0.0238 & 0.0276 & 0.0317\\
7  && 0.0194 & 0.0184 & 0.0170 & 0.0182 & 0.0210 & 0.0237\\
8  && 0.0169 & 0.0153 & 0.0149 & 0.0149 & 0.0172 & 0.0193\\
9  && 0.0149 & 0.0135 & 0.0126 & 0.0137 & 0.0148 & 0.0178\\
10 && 0.0132 & 0.0123 & 0.0116 & 0.0114 & 0.0137 & 0.0149
\end{tabular}
\end{center}
\caption{The mean out of 100 simulations of $\tilde{\mathcal T}_1(\Pi_n,\Pi_{\alpha,\varepsilon})$ for $n=100$.}
\label{fig:mean100}

\end{table}

\begin{table}
\begin{center}
\begin{tabular}{rr|cccccc}
&$ \varepsilon$   & 0.5    & 0.6    & 0.7    & 0.8    & 0.9    & 1.0   \\[-3mm]
$\alpha$ &	&	&	&	&	&	&	\\
\hline
1  && 0.0322 & 0.0367 & 0.0425 & 0.0545 & 0.0743 & 0.0744\\
2  && 0.0147 & 0.0171 & 0.0327 & 0.0275 & 0.0402 & 0.0408\\
3  && 0.0075 & 0.0098 & 0.0189 & 0.0183 & 0.0231 & 0.0263\\
4  && 0.0041 & 0.0062 & 0.0094 & 0.0107 & 0.0174 & 0.0185\\
5  && 0.0031 & 0.0061 & 0.0064 & 0.0089 & 0.0126 & 0.0136\\
6  && 0.0030 & 0.0038 & 0.0057 & 0.0077 & 0.0099 & 0.0115\\
7  && 0.0020 & 0.0035 & 0.0036 & 0.0053 & 0.0070 & 0.0089\\
8  && 0.0016 & 0.0026 & 0.0037 & 0.0040 & 0.0058 & 0.0064\\
9  && 0.0016 & 0.0026 & 0.0040 & 0.0036 & 0.0052 & 0.0057\\
10 && 0.0012 & 0.0025 & 0.0025 & 0.0034 & 0.0041 & 0.0056
\end{tabular}
\end{center}
\caption{The standard deviation out of 100 simulations of $\tilde{\mathcal T}_1(\Pi_n,\Pi_{\alpha,\varepsilon})$ for $n=100$.}
\label{fig:sd100}
\end{table}

\subsection{Climatic time series}

\begin{figure}[p]
\begin{center}
\subfloat[The distribution of the data]{\includegraphics[height=17em,page=2]{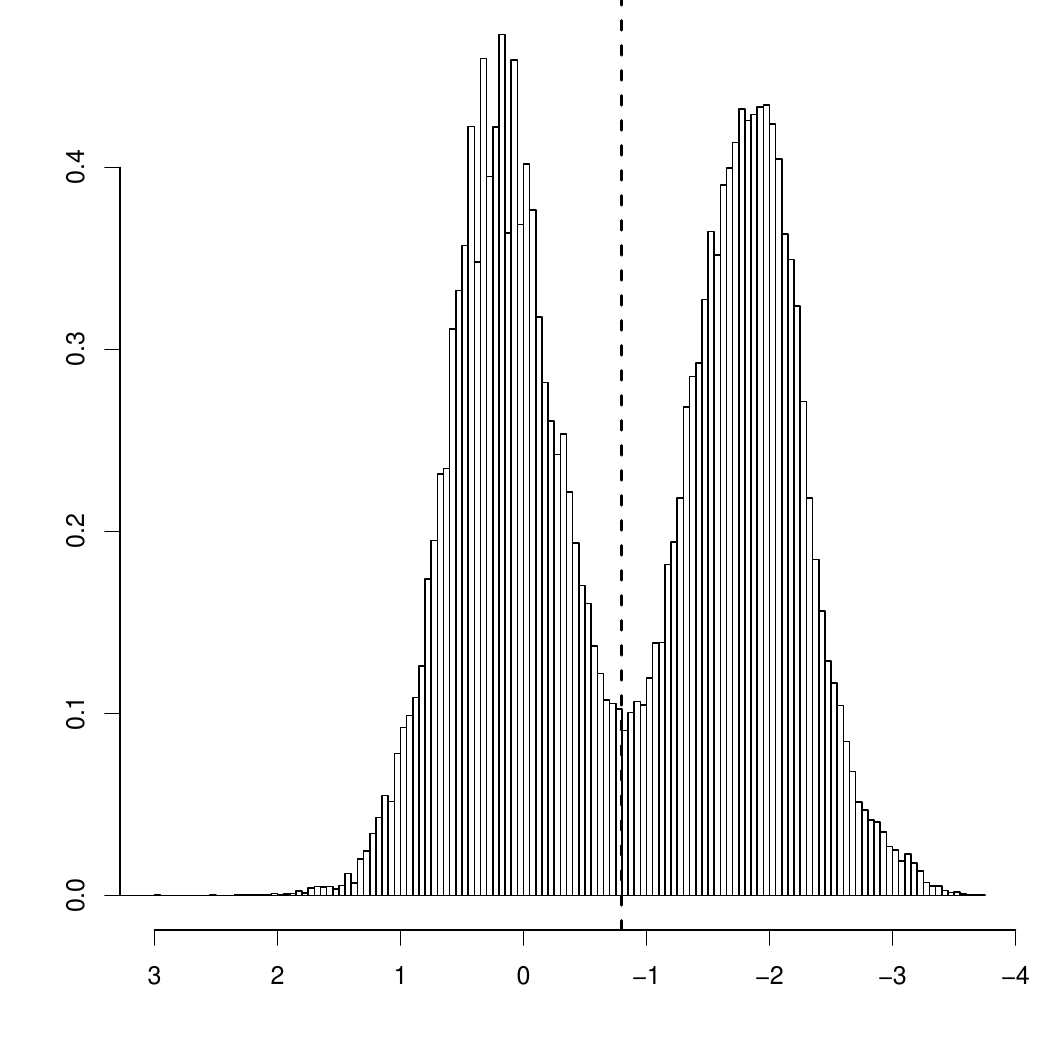}}
\subfloat[Superposition of increments in the cold and warm regime ]{\includegraphics[height=17em,page=3]{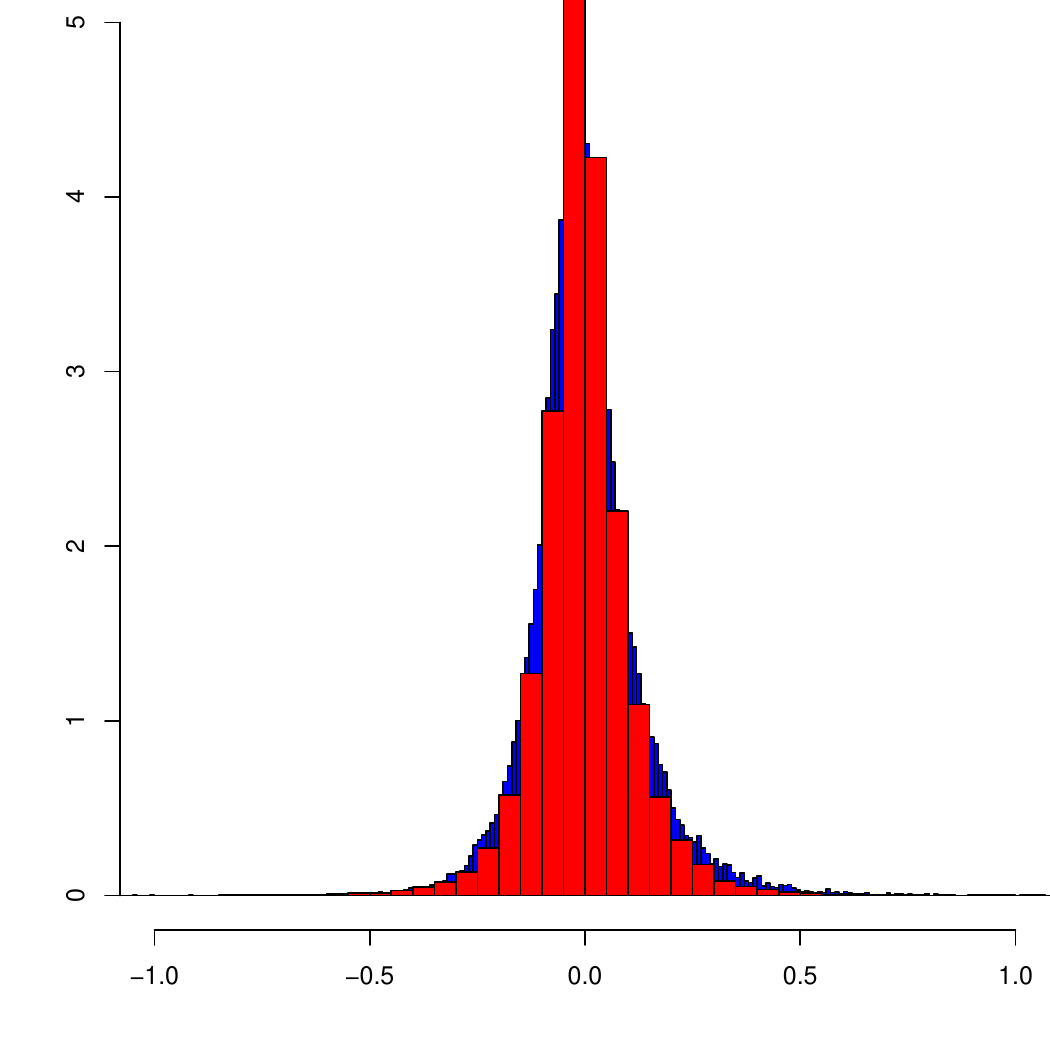}}
\end{center}
\caption{Histograms}
\label{fig:hist}
\end{figure}

Let us turn to the model selection problem for the paleoclimatic time series.
We may extend the preceeding considerations to the dynamical set up of a jump diffusion model.
Figure \ref{fig:hist}(a) shows the histogram of the data and indicates the presence of two well defined regimes, separated by a threshold $s^*\approx -0.8$  centered around the values $0.5$, which corresponds to a \emph{warm} regime, and $-2$, a \emph{cold} regime.
It is therefore reasonable to discretize the state space into these two regimes, separated by a threshold. We assume the characteristics of the noise to be constant in each of the regimes.
Figure \ref{fig:hist}(b) shows the histograms of the increments in each regime.
In both regimes the distributions seem nearly symmetric and may exhibit polynomial tails.

\noindent We will consider a finite intensity version of the model \eqref{eq_ij} given as \eqref{sde3}.
In particular we drop the Gaussian component ($b=0$).
The random measure $\nu$ has a L\'evy kernel $\Pi(x,dv)$ given by a space dependent version of $\Pi_{\alpha,\varepsilon}$ of  Example \ref{ex:powerlaw}, namely
\begin{equation}
\label{eq:kernel1}
\Pi'(x,dv) = \ind(v>\varepsilon^+)\Pi_{\alpha^+(x),\varepsilon^+}(dv) + \ind(v<-\varepsilon^-)\Pi_{\alpha^-(x),\varepsilon^-}(dv)\ .
\end{equation}
The parameters $\alpha^{\pm}$ are assumed to be constant in each (slightly reduced) regime, $x>s^*+\delta, x<s^*-\delta$.
The separation $\delta>0$ is introduced in view of Proposition \ref{prop:1} and Theorem \ref{thm:1} 
to allow for a Lipschitz interpolation of the parameters $\alpha^\pm$,
\begin{equation}
\alpha^\pm(x) =
\begin{cases}
\alpha^{1,\pm} , & x>s^*+\delta,\\
\alpha^{-1,\pm}, & x<s^*-\delta,\\
\text{Lipschitz interpolation}, & 	 \text{otherwise}.
\end{cases}
\label{eq:regimes}
\end{equation}
The class of models is then given by solutions to the SDE
\begin{align}
\label{eq:sde3}
dX(t)
&= a(X(t)) dt + \int_\RR c(X(t-), v) \nu_0(dt, dv)\ ,
\end{align}
where the Cauchy kernel $c$ is constant in each of the regimes \eqref{eq:regimes} and determined by Example \ref{ex:powerlaw}.2.
Between the regimes it inherits the continuous interpolation of the parameters $\alpha^\pm$.
We will now interpret big increments (larger than the threshols $\varepsilon^+$ or smaller 
than $-\varepsilon^-$) of the time series as jumps distributed according to $\Pi(x,\cdot)$.
We extract from the series four sub-samples containing big positive or negative increments in each of the two regimes.
The remaining increments are considered to be continuous and will be neglected.
Theses four sub-samples are now processed in the same way as our simulation study in the previous subsection. 
Their respective empirical measures are compared to a family of L\'evy measures $\Pi_{\alpha^\pm, \e^\pm}$ as before.

\begin{figure}[p]
\begin{center}
\subfloat[Positive jumps in the warm regime:\newline $\alpha_\text{min}=2.8 ,~ \tilde{\mathcal T}_1(\Pi_n,\Pi_{\alpha_\text{min}})=0.048,~n=301$ ]{\includegraphics[height=17em,page=1]{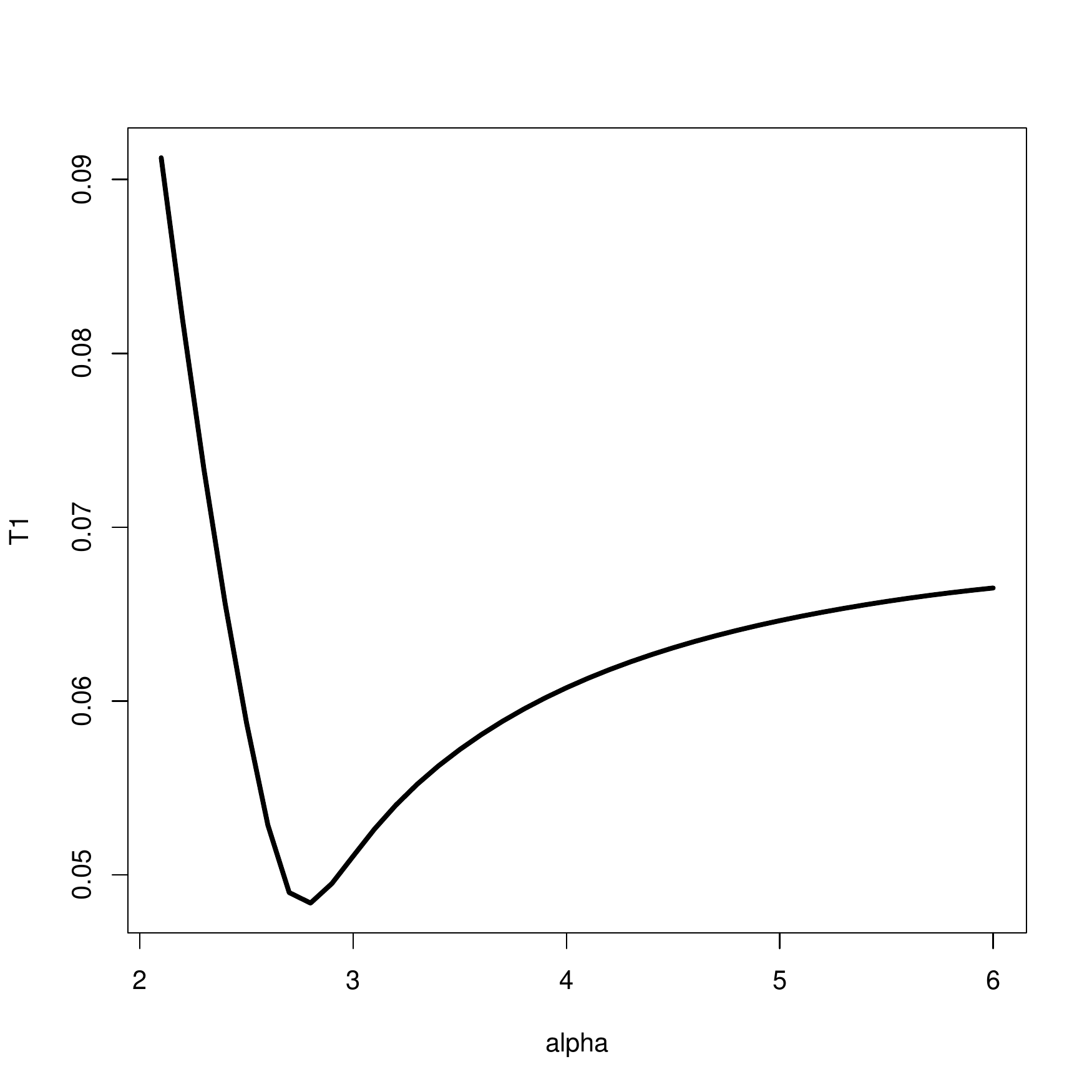}}
\subfloat[Negative jumps in the warm regime:\newline $\alpha_\text{min}=2.9 ,~ \tilde{\mathcal T}_1(\Pi_n,\Pi_{\alpha_\text{min}})=0.050,~n=302$ ]{\includegraphics[height=17em,page=2]{TS2}}\\
\subfloat[Posititve jumps in the cold regime:\newline $\alpha_\text{min}=3.6 ,~ \tilde{\mathcal T}_1(\Pi_n,\Pi_{\alpha_\text{min}})=0.039,~n=593$ ]{\includegraphics[height=17em,page=3]{TS2}}
\subfloat[Negative jumps in the cold regime:\newline $\alpha_\text{min}=4.3 ,~ \tilde{\mathcal T}_1(\Pi_n,\Pi_{\alpha_\text{min}})=0.033,~n=228$ ]{\includegraphics[height=17em,page=4]{TS2}}
\end{center}
\caption{Comparison of the empirical L\'evy measure in the different regimes and different tails to the family $(\Pi_\alpha)_{\alpha\in(2,6)}$.}
\label{fig:AS}
\end{figure}

The thresholds $\varepsilon^\pm$ are chosen in accordance to the polynomial decay of the empirical distribution 
and taken from the preceding study \cite{GHKK15} being $\varepsilon^+=0.36, \varepsilon^-=0.34$ for comparability.
The result is presented in Fig. \ref{fig:AS}.
As can be seen, in each of the cases there is a clear minimizing exponent for the transportation distances varying in $\alpha$.
We find the minimizers
\begin{align*}
&\alpha^{1,+} = 2.8\ ,
&\alpha^{1,-} = 2.9\ , \\
&\alpha^{-1,+}= 3.6\ ,
&\alpha^{-1,-}= 4.3\ .
\end{align*}
We stress that the minimum distances we find are small and of order $0.05$ (note that we consider the normalized 
distance $\tilde{\mathcal T}_1$ with values in $[0,1]$). 
These findings allow to select a simple jump diffusion model 
of type \eqref{eq:sde3} with  L\'evy kernel \eqref{eq:kernel1}.

Theorem \ref{thm:1} allows to control the behavior of the law of this simple model (on path space).
We assume that the drift $a$ in equation \eqref{eq:sde3} is Lipschitz and recall that the interpolation in \eqref{eq:regimes} is Lipschitz, too.
The Lipschitz condition (\ref{Lip0d}) in Theorem \ref{thm:1} for $\Pi'$ follows from Proposition~\ref{prop:1}.
Theorem \ref{thm:1} guarantees that the laws of such jump diffusions on a reasonable time horizon are controlled by the transportation distance of the L\'evy measures.
Hence the jump diffusions selected by the 
values obtained by the minimization procedure provide a good estimate 
on the space of such models. 

It is possible to obtain the rate of convergence of the empirical L\'evy measure $\Pi_n$ under the transportation 
distance similarly to the case of coupling distances (c.f. \cite{GHKK15}).
In this study we content ourselves to the convincing simulation results of Section \ref{sec:simulation}.
As already mentioned a similar fitting approach for a family of polynomial L\'evy measures to this time-series is adopted in \cite{GHKK15}. 
The frame work there is restricted to models with additive noise and does not allow for spacial dependence of the jump kernel 
and cannot distinguish different regimes.
However, the polynomial decay of the left and right tail of the jump measure is studied.
As a result the minimizing exponents for the coupling-distance there are found 
to be $\alpha^+= 3.6$ (positive jumps) and $\alpha^-=3.55$ (negative jumps).

\noindent Our findings seem to be supportive with this result given that we consider different spacial regimes.
Combining the negative jumps in the different regimes should lead to an average of the exponents $\alpha^{1,-},\alpha^{-1,-}$.
If we average the exponents weighted by the number of negative jumps in the corresponding regime we obtain
\begin{equation*}
\frac{2.9\times 302 + 4.3 \times 228}{302+228} = 3.5 \approx 3.55 \ .
\end{equation*}
For the positive jumps averaging behaves slightly worse and we obtain roughly
\begin{equation*}
\frac{2.8\times 301 + 3.6 \times 593}{301+593} = 3.33 \approx 3.6 \ .
\end{equation*}
In any case this study confirms earlier findings in \cite{GHKK15} of a robust polynomial tail behavior with exponents $\alpha$ well beyond $2$.
Originally a jump diffusion model with $\alpha$-stable noise component was proposed in \cite{Dit99a} which restricts $\alpha$ to values in $(0,2)$.
A waiting time analysis of 42 transition between warm and cold temperature regimes pointed to $\alpha\approx 1.75$.
Later a similar model was investigated by refined methods in \cite{GI14} that supported an exponent of $\alpha\approx 1.75$ within the framework of $\alpha$-stable L\'evy noise.
We stress that the class we are considering here are generic in the class of heavy tailed L\'evy measures and approximate $\alpha$-stable ones.
The findings of $\alpha$ being close to 2 within the class of $\alpha$-stable distributions 
may also indicate a possible exponent beyond 2, since their methods perform poorly for $\alpha$ approaching 2.

\section{Proofs}\label{sec: proofs}

\begin{proof} \emph{of Theorem \ref{est_forT_2_rho}: }\\
The strategy of the proof is to bound first 
$\EE\big[\left(X_{1}(t)-X_{2}(t)\right)^{2}\wedge 1\big]$ 
with the help of the cut-off technique developed in
\cite{GHKK14}. Using the elementary inequality
\be\label{arctan}
y\wedge 1\lqq {4\over \pi}\arctan y, \quad \mbox{ for any }y\gqq 0,
\ee
we obtain for $F(y) = \arctan(y^2)$
\[
\rho^2(x,y) \lqq {4\over \pi} F(|x-y|), \qquad \mbox{ for all } x,y \in \RR.
\]
Applying It\^o's formula and Doob's maximal inequality to
$
F(|X_1-X_2|)
$
we establish Gronwall estimates for the process $F(|X_1-X_2|).$
As a preparation for the two dimensional It\^o formula for $F(|X_1-X_2|)$
we rewrite the process $Y(t)=X_{1}(t)-X_{2}(t)$ 
and 
$a^*_i(x) = a_i(x) + \bar a_i(x)$ 
with common Lipschitz constant $L_{a^*} = L_{a} + L_{\bar a}$ as follows
\begin{align*}
d(X_1-X_2)(t) =&
(a^*_1(X_1) - a^*_2(X_2)) dt + (b_1(X_1) - b_2(X_2)) dW(t)\\
&+ \int_{|v|> 1} (c_1(X_1(t-),v) - c_2(X_2(t-),v) \nu_0(dv,dt)\\
&- \int_{|v|>1} \Big(\tau(c_1(X_1(t-),v)) - \tau(c_2(X_2(t-),v))\Big)\Pi_0(dv)dt\\
&-  \int_{|v|\lqq 1}  \Big(\tau(c_1(X_1(t-), v)) - c_1(X_1(t-), v)\Big) \Pi_0(dv) dt \\
&+  \int_{|v|\lqq 1}  \Big(\tau(c_2(X_2(t-), v)) - c_2(X_2(t-), v)\Big) \Pi_0(dv) dt\\
&-\int_{|v|\lqq 1} c_1(X_1(t-),v)  \ti \nu_0(dv,dt) + \int_{|v|\lqq 1} c_2(X_2(t-),v)  \ti \nu_0(dv,dt),
\end{align*}
where
$$\tau(v)=(|v|\wedge 1)\,\mathrm{sign}\,(v), \quad  v\in \RR.$$
Applying It\^o's formula given for instance in Chapter 2 in \cite{Vat-Ik} we get
\begin{align*}
d F(Y(t)) &= F'(Y(t-)) (a_1^*(X_1) - a_2^*(X_2)) dt + F'(Y(t))(b_1(X_1)-b_2(X_2)) dW(t)\\
&\qquad+ \frac{1}{2} F''(Y(t))(b_1(X_1)-b_2(X_2))^2 dt   \\
&\qquad-\int_{|v|\lqq 1}  F'(Y(t-)) \Big(\tau(c_1(X_1(t-), v)) - c_1(X_1(t-), v)\Big) \Pi_0(dv) dt \\
&\qquad+\int_{|v|\lqq 1}  F'(Y(t-))\Big(\tau(c_2(X_2(t-), v)) - c_2(X_2(t-), v)\Big) \Pi_0(dv) dt\\
&\qquad- \int_{|v|>1} \Big(F'(Y(t-))(\tau(c_1(X_1(t-),v)) - \tau(c_2(X_2(t-),v)))\Big)\Pi_0(dv)dt\\
&\qquad+ \int_{|v|>1} \Big(F(Y(t-)+(c_1(X_1(t),v)-c_2(X_2(t),v))) -F(Y(t-))\Big) \nu_0(dv,dt) \\
&\qquad + \int_{|v|\lqq 1} \Big(F(Y(t-)+(c_1(X_1(t),v)-c_2(X_2(t),v))) -F(Y(t-))\Big) \ti \nu_0(dv,dt) \\
&\qquad + \int_{|v|\lqq 1} \Big(F(Y(t-)+(c_1(X_1(t),v)-c_2(X_2(t),v))) -F(Y(t-)) \\
&\qquad\qquad- F'(Y(t-))(c_1(X_1(t),v)-c_2(X_2(t),v))\Big) \Pi_0(dv) dt.
\end{align*}
After rearrangements, we finally obtain the representation
\ba\label{F}
F(Y(t))=F(|x_1- x_2|) + M(t)+\int_{0}^{t}\int_{\RR} g(Y(s),X_{1}(s),X_{2}(s))\Pi_{0}(dv)ds,
\ea
where
\ba
M(t)=&\int_{0}^{t}F'(Y(s))(b_1(X_1(s))-b_2(X_2(s))) dW(s)\\
&+\int_{0}^{t}\int_{\RR} \Big(F(Y(s-)+(c_1(X_1(s),v)-c_2(X_2(s),v))) -F(Y(s-)) \Big)[\nu_0(dv,ds)-\Pi_{0}(dv)ds]
\\
\ea
and
\ba\label{g}
g(y,u_{1},u_{2})&=F'(y) (a_1^*(u_1) - a_2^*(u_2))+\frac{1}{2} F''(y)(b_1(u_1)-b_2(u_2))^2\\
&+\int_{\RR} [F(y+(c_1(u_{1},v)-c_2(u_{2},v))) -F(y)- F'(y)(\tau(c_1(u_{1},v)) - \tau(c_2(u_{2},v)))]\Pi_0(dv)\\
&=:g_{1}(y,u_{1},u_{2})+g_{2}(y,u_{1},u_{2})+g_{3}(y,u_{1},u_{2}),\quad (y,u_{1},u_{2})\in\RR^{3},\quad y=u_{1}-u_{2}.
\ea
For the function $F$ and its derivatives, we have  the following bounds taken from \cite{GHKK14}
\begin{align}
&F'(y)={2y\over 1+y^4}, \quad \mbox{ implying }\qquad |F'(y)|\lqq \frac{3^{3/4}}{2} \mbox{ for all }y \in \RR,\label{dF}\\ 
&F'(y)y={2y^2\over 1+y^4}\lqq (2y^2)\wedge 1 = 2 (y^2\wedge \frac{1}{2})\lqq \frac{F(y)}{\arctan(1/2)},\label{dFy}\\
&F''(y)=2{1-3y^4\over (1+y^4)^2}, \quad \mbox{ implying } \qquad |F''(y)|\lqq 2, \label{ddF}\\
&|F(y+\delta)-F(y)-F'(y)\delta|\lqq {\delta^2\over 2}\sup_v|F''(v)|\lqq \delta^2.\label{F-F-dF}
\end{align}
Now we can bound each summand of (\ref{g}).
Using (\ref{Lip0a}) and (\ref{dFy}) we obtain
\ba\label{a_estimate}
g_{1}(y,u_{1},u_{2})\lqq \frac{L_{a^*}}{\arctan(1/2)}F(y) +\frac{3^{3/4}}{2} \|a_1^* - a_2^*\|_\infty 
\qquad \forall u_1, u_2, y = u_1-u_2 \in \RR.
\ea
Then by (\ref{Lip0b}) $g_{2}$ is not greater than
\ba\label{b_estimate}
g_{2}(y,u_{1},u_{2})\lqq (b_1(u_1)-b_2(u_2))^2 \qquad \forall u_1, u_2, y = u_1-u_2 \in \RR.
\ea
The estimate of $g_{3}$ is slightly more involved
and for convenience we introduce for fixed $u_i$ the notation
\begin{align*}
\zeta_i(v) := c_i(u_{i},v), \qquad v\in \RR.
\end{align*}
Let us first consider the integral $g_3$ on the set $\{v \in \RR~|~|\zeta_1(v) - \zeta_2(v)|\lqq 1\}$
and apply
\begin{align*}
F(y + (z_1-z_2)) - F(y) \lqq F'(y)(z_1-z_2) + (z_1-z_2)^2
\end{align*}
to the first summand of
\begin{align*}
&\int_{\{v \in \RR~|~|\zeta_1(v) - \zeta_2(v)|\lqq 1\}}
\Big([F(y+(\zeta_1(v) - \zeta_2(v))) -F(y)] - F'(y)[\tau(\zeta_1(v)) - \tau(\zeta_2(v))]\Big) \Pi_0(dv)\\[3mm]
&\lqq \int_{\{v \in \RR~|~|\zeta_1(v) - \zeta_2(v)|\lqq 1\}}
\Big(F'(y)\Big[ \zeta_1(v) - \zeta_2(v)- \tau(\zeta_1(v)) + \tau(\zeta_2(v))\Big]
+ (\zeta_1(v)-\zeta_2(v))^2\Big) \Pi_0(dv).
\end{align*}
Observe that the term in square brackets vanishes if both, $|\zeta_i(v)|\lqq 1$ for $i=1,2$.
Hence the first integral reads
\begin{align*}
&\int_{\{v \in \RR~|~|\zeta_1(v) - \zeta_2(v)|\lqq 1\}}
F'(y)\Big[ \zeta_1(v) - \zeta_2(v)- \tau(\zeta_1(v)) + \tau(\zeta_2(v))\Big] \Pi_0(dv)\\
&= \int_{\substack{
\{v \in \RR ~|~|\zeta_1(v) - \zeta_2(v)|\lqq 1 \\
\quad \mbox{ \tiny and } |\zeta_1(v)|>1 \mbox{ \tiny or } |\zeta_2(v)|>1\}}}
F'(y)\Big[ \zeta_1(v) - \zeta_2(v)- \tau(\zeta_1(v)) + \tau(\zeta_2(v))\Big] \Pi_0(dv)\\
&\lqq 2 F'(y) \int_{\substack{
\{v \in \RR ~|~|\zeta_1(v) - \zeta_2(v)|\lqq 1  \\
\quad \mbox{ \tiny and } |\zeta_1(v)|>1 \mbox{ \tiny or } |\zeta_2(v)|>1\}}}
\Big[ |\zeta_1(v)-\zeta_2(v)|\wedge 1\Big] \Pi_0(dv)\\
&\lqq 2 F'(y) \Big(\int_{\RR}
\Big[ |\zeta_1(v)-\zeta_2(v)|\wedge 1\Big]^2 \Pi_0(dv)\Big)^\frac{1}{2}
\Pi_0\Big(\{v\in \RR~|~|\zeta_1(v)|>1\mbox{ or }|\zeta_2(v)|>1\}\Big)^\frac{1}{2} \\
&\lqq 2 F'(y)  \mathcal T_2 (\Pi_1(u_{1}, \cdot), \Pi_2(u_{2}, \cdot))
\Big(\Pi_0\big(\{v\in \RR~|~|\zeta_1(v)|>1\}\big)+ \Pi_0\big(\{v\in \RR~|~|\zeta_2(v)|>1\}\big)\Big)^\frac{1}{2},
\end{align*}
where we have used that $|\zeta_1(v) - \zeta_2(v)|\lqq 1$ implies
\begin{align*}
|\zeta_1(v)- \zeta_2(v)| + |\tau(\zeta_1(v)) - \tau(\zeta_2(v))|\lqq 2 (|\zeta_1(v) - \zeta_2(v)|\wedge 1).
\end{align*}
We consider now the integral $g_3$ on the remainder set $\{v \in \RR~|~|\zeta_1(v) - \zeta_2(v)|> 1\}$.
By its definition in (\ref{g})
$g_3$ is uniformly bounded by $\pi + \frac{3^\frac{3}{4}}{2}$
due to the boundedness of $F$ and (\ref{dF}).
Therefore the function $g_{3}$ in (\ref{g}) is bounded by the sum
\begin{align}
g_{3}(y,u_{1},u_{2})&\lqq 2 F'(y)  \mathcal T_2(\Pi_1(u_{1}, \cdot), \Pi_2(u_{2}, \cdot))
\Big(\Pi_0(\{v:|\zeta_1(v)|>1\}+ \Pi_0(\{v:|\zeta_2(v)|>1\})\Big)^\frac{1}{2}\nonumber\\
&\qquad + (\pi + \frac{3^\frac{3}{4}}{2}) \int_{\RR} ((\zeta_1(v) - \zeta_2(v))^2 \wedge 1) \Pi_0(dv) \nonumber\\
&\lqq 2 F'(y)  \mathcal T_2(\Pi_1(u_{1}, \cdot), \Pi_2(u_{2}, \cdot)) \Big(\Pi_0(\{v:|\zeta_1(v)|>1\}
+ \Pi_0(\{v:|\zeta_2(v)|>1\})\Big)^\frac{1}{2} \nonumber\\
&\qquad + (\pi + \frac{3^\frac{3}{4}}{2}) \mathcal T_2(\Pi_1(u_{1}, \cdot), \Pi_2(u_{2}, \cdot))^2.\label{g_3}
\end{align}
Summarising (\ref{a_estimate}), (\ref{b_estimate}) and (\ref{g_3}) we get for the process $F(Y(t))$ the following estimate almost surely
\begin{align*}
&F(Y(t))\\
& \lqq \frac{\pi}{2} \rho^2(x_1, x_2) + T \frac{3^{3/4}}{2} \|a_1^*-a_2^*\|_\infty -\frac{L_{a^*}}{\arctan(1/2)} \int_0^tF(Y(s)) ds +\int_0^t \Big((b_1(X_{1}(s)) - b_2(X_{2}(s)))^2 ds\\
&\qquad + 2  \int_0^t F'(Y(s-))\mathcal T_2(\Pi_1(x_{1}, \cdot), \Pi_2(x_{2}, \cdot))
\Big(\Pi_0(\{v:|\zeta_1(v)|>1)^\frac{1}{2}+ \Pi_0(\{v:|\zeta_2(v)|>1\})^\frac{1}{2}\Big) ds\\
&\qquad + \int_0^t \Big(\pi + \frac{3^\frac{3}{4}}{2})\mathcal T_2^{2}(\Pi_1(X_{1}(s-), \cdot), \Pi_2(X_{2}(s-), \cdot)\Big) ds + M(t).
\end{align*}
In the sequel we apply the Lipschitz continuity of the coefficients $b_{i},\ i=1,2$ and the kernels $\Pi_i,\ i=1,2$, 
however taken at different points.
By (\ref{Lip0b}) and (\ref{Lip0c}) we obtain $\om$-wise
\begin{align*}
&\mathcal T_2(\Pi_1(X_{1}(s-), \cdot), \Pi_2(X_{2}(s-), \cdot))\lqq L_{\Pi} (|Y(s-)|\wedge 1) + \sup_{x\in\RR}\mathcal T_2(\Pi_1(x, \cdot), \Pi_2(x, \cdot))
\end{align*}
and by (\ref{arctan})
\begin{align}
&\mathcal T_2^{2}(\Pi_1(X_{1}(s-), \cdot), \Pi_2(X_{2}(s-), \cdot))\lqq \frac{8L_{\Pi}}{\pi} F(Y(s-)) + 2\sup_{x\in\RR}\mathcal T_2^{2}(\Pi_1(x, \cdot), \Pi_2(x, \cdot)),
\nonumber\\
&(b_1(X_{1}) - b_2(X_{2}))^{2}\lqq {8L_{b}\over \pi} F(Y(s))+2\sup_{x\in\RR}(b_1(x) - b_2(x))^{2}.\label{T,b}
\end{align}
Setting
\[
K_1 := \sup_{u_1, u_2 \in \RR} \Big(\Pi_0(\{v:|c_1(u_1, v)|>1)+ \Pi_0(\{v:|c_2(u_2, v)|>1\})\Big)^\frac{1}{2}
\]
and the help of (\ref{dFy}) we obtain almost surely 
\begin{align}
F(Y(t))
&\lqq \frac{\pi}{2} \rho^2(x_1, x_2)+ T \frac{3^{3/4}}{2} \|a_1^*-a_2^*\|_\infty + \frac{L_{a^*}}{\arctan(1/2)} \int_0^t F(Y(s)) ds + {8L_{b}\over \pi}\int_0^t F(Y(s)) ds\nonumber\\
&\qquad + 2t\sup_{x\in\RR}(b_1(x) - b_2(x))^{2} \nonumber\\
&\qquad +M(t) + 2  K_{1}\int_0^t |F'(Y(s))| |Y(s)|  ds + 4K_{1}t\sup_{x\in\RR}\mathcal T_2(\Pi_1(x, \cdot), \Pi_2(x, \cdot))\nonumber\\
&\qquad + (8 + \frac{4\cdot3^\frac{3}{4}}{\pi})L_{\Pi}\int_0^t F(Y(s)) ds
+ (2\pi + 3^\frac{3}{4})t\sup_{x\in\RR}\mathcal T_2^{2}(\Pi_1(x, \cdot), \Pi_2(x, \cdot))\nonumber\\
&\lqq \frac{\pi}{2} \rho^2(x_1, x_2) + \Big(\frac{L_{a^*}+2K_{1}}{\arctan(1/2)} + {8L_{b}\over \pi}
+ (8 + \frac{4\cdot3^\frac{3}{4}}{\pi})L_{\Pi}\Big)\int_0^t F(Y(s))ds+ M(t) \nonumber\\
&\qquad + 2t\sup_{x\in\RR}(b_1(x) - b_2(x))^{2}
+ 4tK_{1}\sup_{x\in\RR}\mathcal T_2(\Pi_1(x, \cdot), \Pi_2(x, \cdot))\nonumber\\
&\qquad + T \frac{3^{3/4}}{2} \|a_1^*-a_2^*\|_\infty +t(2\pi + 3^\frac{3}{4})\sup_{x\in\RR}\mathcal T_2^{2}(\Pi_1(x, \cdot), \Pi_2(x, \cdot)).\label{pre_Gronwall}
\end{align}
We define 
$$
Q:=\frac{L_{a^*}+2K_{1}}{\arctan(1/2)} + {8L_{b}\over \pi} + (8 + \frac{4\cdot3^\frac{3}{4}}{\pi})L_{\Pi},
$$
take the expectation and apply Gronwall's lemma for $t\in [0, T]$
\ba\label{GronwallF}
\EE[F(Y(t))]
&\lqq \Big(\frac{\pi}{2} \rho^2(x_1, x_2) + T \frac{3^{3/4}}{2} \|a_1^*-a_2^*\|_\infty + 2T\sup_{x\in\RR}(b_1(x) - b_2(x))^{2} \\
&\qquad+ 4TK_{1}\sup_{x\in\RR}\mathcal T_2(\Pi_1(x, \cdot), \Pi_2(x, \cdot)) + T (2\pi + 3^\frac{3}{4})\sup_{x\in\RR}\mathcal T_2^{2}(\Pi_1(x, \cdot), \Pi_2(x, \cdot))\Big)e^{QT}.
\ea
Now establish the estimate for $\EE[\sup_{t \in [0, T]} F(Y(t))]$.
For that purpose we take the supremum over both sides of (\ref{pre_Gronwall})
and observe that
\begin{align}
&\EE[\sup_{t \in [0, T]} F(Y(t))] 
\lqq \frac{\pi}{2} \rho^2(x_1, x_2) + Q\int_0^T \EE[ F(Y(s))]ds + \EE[\sup_{t\in [0, T]} |M(t)|] +\nonumber\\
& \qquad T\Big(\frac{3^{3/4}}{2} \|a_1^*-a_2^*\|_\infty +2\sup_{x\in\RR}(b_1(x) - b_2(x))^{2} 
+ 4K_{1}\sup_{x\in\RR}\mathcal T_2(\Pi_1(x, \cdot), \Pi_2(x, \cdot))\nonumber\\
&\qquad \quad +(2\pi 
+ 3^\frac{3}{4})\sup_{x\in\RR}\mathcal T_2^{2}(\Pi_1(x, \cdot), \Pi_2(x, \cdot))\Big).\label{presup}
\end{align}
For the integrand of the first integral we apply (\ref{GronwallF}).
In order to bound the martingal term we use Doob's maximum inequality 
\begin{align*}
\EE[\sup_{t\in [0, T]} |M(t)|]
&\lqq \EE[\sup_{t\in [0, T]} |M(t)|^2]^\frac{1}{2}
\lqq 4 \EE[|M(T)|^2]^\frac{1}{2}
\end{align*}
and rewrite $\EE[|M(T)|^2]$ as follows
\begin{align*}
\EE[|M(T)|^2]
&=\EE\Big[\Big(\int_0^T\int_\RR \Big[F\big(Y(t-)+(c_1(X_1(t),v)-c_2(X_2(t),v))\big) -F(Y(t-))\Big] \nu_0(dv,dt) \\
&\qquad -\int_0^T \int_\RR [F\big(Y(t-)+(c_1(X_1(t),v)-c_2(X_2(t),v))\big) -F(Y(t-))] \Pi_0(dv)dt\Big)\Big)^2\Big]\\
&\qquad + \int_0^T \EE\bigg[\Big( F'(Y(t))\Big(b_1(X_1(t)) - b_2(X_2(t))\Big)\Big)^2\bigg] dt \nonumber\\
&= \int_0^T \int_\RR \EE\bigg[\Big(F\big(Y(t-)+(c_1(X_1(t),v)-c_2(X_2(t),v))\big) -F(Y(t-))\Big)^2\bigg] \Pi_0(dv)dt\\
&\qquad + \int_0^T \EE\bigg[\Big( F'(Y(t))\Big(b_1(X_1(t)) - b_2(X_2(t))\Big)\Big)^2\bigg] dt.
\end{align*}
The analogous separation argument as for the term $g_3$
of the cases $\{v\in \RR~|~|c_1(x, v) - c_2(x, v)| \lqq 1\}$ and its complement
and using $F(y+\delta) - F(y) \lqq 2 \delta$ for $\delta\in (0,1)$ and $y\in \RR$
we obtain with the help of (\ref{Lip0b}), (\ref{Lip0c}) and (\ref{arctan}) the estimate
\begin{align*}
\EE[|M(T)|^2]
&\lqq 2\int_0^T \int_\RR \EE\Big[(c_1(X_1(t),v)-c_2(X_2(t),v))^2\wedge 1\Big] \Pi_0(dv)dt\\
&\quad +2 \int_0^T \EE\Big[\Big(b_1(X_1(t)) - b_2(X_2(t))\Big)^2\Big] dt\\
&= 2\int_0^T \EE\Big[ \mathcal T_2^{2}(\Pi_1(X_1(t), \cdot), \Pi_2(X_2(t), \cdot)\Big] dt
+2\int_0^T \EE\Big[\Big(b_1(X_1(t)) - b_2(X_2(t))\Big)^2\Big] dt\\
&\lqq \frac{16L_{\Pi}}{\pi}\int_0^T \EE\Big[ F(Y(t))\Big] dt+4T\sup_{x\in\RR}\mathcal T_2^{2}(\Pi_1(x, \cdot), \Pi_2(x, \cdot))\\
&\quad +\frac{16L_{b}}{\pi}\int_0^T \EE\Big[ F(Y(t))\Big] dt+4T\sup_{x\in\RR}(b_1(x) - b_2(x))^{2}.
\end{align*}
The last inequality comes from the triangle inequality (\ref{T,b}).
Finally inserting (\ref{GronwallF}) we obtain
\begin{align*}
&\EE\Big[\sup_{t \in [0, T]} F(Y(t))\Big] \\
&\lqq (T^{2} Qe^{QT}+T)\Bigg(\frac{\pi}{2} \rho^2(x_1, x_2) + \frac{3^{3/4}}{2} \|a_1^*-a_2^*\|_\infty +2\sup_{x\in\RR}(b_1(x) - b_2(x))^{2} \\
&\qquad + 4K_{1}\sup_{x\in\RR}\mathcal T_2(\Pi_1(x, \cdot), \Pi_2(x, \cdot))
+(2\pi + 3^\frac{3}{4})\sup_{x\in\RR}\mathcal T_2^{2}(\Pi_1(x, \cdot), \Pi_2(x, \cdot))\Bigg)\\
& \qquad +\Bigg( 32 (L_\Pi + L_b)  \rho^2(x_1, x_2) + \frac{64 (L_\Pi + L_b)}{2\pi} T \|a_1^*-a_2^*\|_\infty^2 \\
&\qquad \quad + \left(\frac{32}{\pi} (L_{\Pi}+L_{b}) T^{2} e^{QT}+4T\right)\sup_{x\in\RR}(b_1(x) - b_2(x))^{2} \\
&\qquad \quad + \frac{64}{\pi} (L_{\Pi}+L_{b}) K_{1}e^{QT}\sup_{x\in\RR}\mathcal T_2(\Pi_1(x, \cdot), \Pi_2(x, \cdot))\\
&\qquad \quad +\Big(\frac{16}{\pi} (L_{\Pi}+L_{b}) (2\pi + 3^\frac{3}{4})T^{2}e^{QT}+4T\Big)\sup_{x\in\RR}\mathcal T_2^{2}(\Pi_1(x, \cdot), \Pi_2(x, \cdot))\Bigg)^{1/2}.
\end{align*}
Taking the maximum of all the constants which we denote by $D = D(K_{1},L_{a^*},L_{b},L_{\Pi},T)$, 
using $x + \sqrt{x} \lqq 2 G(x)$ for any $x\in \RR$, 
and rearranging the terms we obtain the desired estimate (\ref{thm_estim}). This finishes the proof.
\end{proof}

\vspace{2cm}
\begin{proof} \emph{of the Theorem \ref{thm:1}}\\
In order to get estimate (\ref{thm_estim_for_1_rho}) we may use the same scheme as 
in the proof of the Theorem (\ref{est_forT_2_rho}) 
but for a $\cC^1$-function $H^\delta:\RR \ra \RR_+$ instead of $F$ satisfying 
\begin{align}
&\|H^\delta\|_\infty \lqq 2 , \qquad \|(H^\delta)'\|_\infty \lqq 2, \qquad \lim_{\delta\ra 0+} \|H^\delta(y) - (1\wedge |y|)\|_\infty = 0,
\end{align}
since there is no Brownian part and the intensity of the jump part is finite. 
We take the process $Y(t)=X_{1}(t)-X_{2}(t),$ which is the difference of the processes defined 
in (\ref{sde3}), 
and apply It\^o's formula for $H^\delta(Y)$. 
Analogously to (\ref{F}) we obtain for $H^\delta(Y)$ the simpler formula
\ba\label{H}
H^\delta(Y(t))&= \int_0^t \Big(a_1(X_1(s)) - a_2(X_2(s))\Big) ds  + M(t)\\
&\qquad +\int_{0}^{t}\int_{\RR} \big[H^\delta(Y(s)+(c_1(X_{1}(s),v)-c_2(X_{2}(s),v)))-H^\delta(Y(s))\big]\Pi_0(dv)ds+H^\delta(Y(0)),
\ea
where the martingale term $M^\delta(t)$ is given as 
$$
M^\delta(t)=\int_{0}^{t}\int_{\RR} \Big[H^\delta(Y(s-)+(c_1(X_1(s),v)-c_2(X_2(s),v))) -H^\delta(Y(s-))\Big][\nu_0(dv,ds)-\Pi_{0}(dv)ds].
$$
We obtain the estimate for the inner integral 
in the second summand of (\ref{H}) 
on two separate sets 
\begin{equation}\label{eq: separation}
\{v: |c_1(X_{1}(s),v)-c_2(X_{2}(s),v)|\lqq 1\} \mbox{ and } \{v: |c_1(X_{1}(s),v)-c_2(X_{2}(s),v)|> 1\}
\end{equation}
for fixed $s\in[0,t]$. Using the fact that
$
|(H^\delta)'(y)|\lqq 2 \label{dH}
$
and a Taylor expansion of $H^\delta$ in $y$ yields 
on the set $\{v\in \RR: |c_1(X_{1}(s),v)-c_2(X_{2}(s),v)|\lqq 1\}$ almost surely 
\begin{align}
|H^\delta(y+(c_1(x_{1},v)-c_2(x_{v},v)))-H^\delta(y)|\lqq 2|c_1(x_{1},v)-c_2(x_{2},v))|\label{Taylor_H}.
\end{align}
Hence we obtain the following $\om$-wise estimates
\begin{align}
H^\delta(Y(t))
&\lqq M^\delta(t)+ 2 L_a \int_0^t (|Y(s)| \wedge 1) ds + 2 t \|a_1 -a_2\|_\infty  \\
&\quad +2\int_{0}^{t}\int_{\{v: |c_1(X_{1}(s),v)-c_2(X_{2}(s),v)|\lqq 1\}}|c_1(X_{1}(s),v)-c_2(X_{2}(s),v))|~\Pi_0(dv)ds\nonumber\\
&\quad +4\int_{0}^{t}\int_{\{v: |c_1(X_{1}(s),v)-c_2(X_{2}(s),v)|> 1\}}~\Pi_0(dv)ds+H^\delta(Y(0))\nonumber\\
&\lqq M^\delta(t)+ 2 L_a \int_0^t (|Y(s)| \wedge 1) ds  + 2 t \|a_1 -a_2\|_\infty \nonumber\\
&\quad +4\int_{0}^{t}\int_{\RR}\big(|c_1(X_{1}(s),v)-c_2(X_{2}(s),v))|\wedge 1\big)~\Pi_0(dv)ds+H^\delta(Y(0))\nonumber\\
&= M^\delta(t)+ 4 L_a \int_0^t (|Y(s)| \wedge 1) ds  + 2 t \|a_1 -a_2\|_\infty \nonumber\\ 
&\quad +4\int_{0}^{t}\mathcal T_1(\Pi_{1}(X_{1}(s),\cdot),\Pi_{2}(X_{2}(s),\cdot)) ds+H^\delta(Y(0)),
\label{estm: pre H}
\end{align}
where we have used the boundedness of the function $H^\delta$ and the fact that $2|y|\wedge 4\lqq 4(|y|\wedge1)$ for $y\in\RR.$
The inner integral in the last expression is by definition $\mathcal T_1(\Pi_{1}(X_{1}(s),\cdot),\Pi_{2}(X_{2}(s),\cdot)).$ 
For the second term on the right-hand side of (\ref{estm: pre H}) we apply 
the triangular inequality for $\mathcal T_1$ and the Lipschitz condition (\ref{Lip0d}). 
We take the expectation on (\ref{estm: pre H}), send $\delta \ra 0$ and obtain 
\begin{align}
\EE \big[|Y(t)|\wedge 1\big]
&\lqq (2L_a + 4 L_{\Pi})\int_{0}^{t}\EE\big[|Y(s)|\wedge 1\big]ds+ 2 t \|a_1 -a_2\|_\infty
+4T\sup_{x\in\RR}\mathcal T_1(\Pi_{1}(x,\cdot),\Pi_{2}(x,\cdot))\nonumber\\
&\qquad +\rho(x_1, x_2).\label{pre Gronwall 1}
\end{align}
Gronwall's lemma for $t\in[0, T]$ and the monotonicity on the right-hand side imply 
\ba\label{gronwall_first_mom}
\EE \left[|Y(t)|\wedge 1\right]\lqq e^{(2 L_a + 4 L_{\Pi})T}\big(2 T \|a_1 -a_2\|_\infty 
+ 4T\sup_{x\in\RR}\mathcal T_1(\Pi_{1}(x,\cdot),\Pi_{2}(x,\cdot))+\rho(x_1, x_2)\big).
\ea
In what follows we denote $\Delta:= 2 T \|a_1 -a_2\|_\infty + 4T\sup_{x\in\RR}\mathcal T_1(\Pi_{1}(x,\cdot),\Pi_{2}(x,\cdot))+\rho(x_1, x_2).$
Further, taking the supremum in $t\in [0, T]$ on both sides of (\ref{estm: pre H}) we 
use the monotonicity on the right-hand side and take the expectation 
\ba
\EE\big[\sup_{t\in[0,T]} H^\delta(Y(s))\big]
\lqq \EE[\sup_{t\in[0,T]}|M^\delta_{t}|]+
(2 L_a + 4 L_{\Pi})\int_{0}^{T}\EE\big[\left(|Y(s)|\wedge 1\right)\big]ds + \Delta.\label{pre sup}
\ea
For the martingale term in (\ref{pre sup}) we apply Doob's maximal moment inequality with $p=2$, 
which yields 
$$
\EE\sup_{t\in [0,T]}|M^\delta_t|^2\lqq 4\EE|M^\delta_T|^2
$$
and hence 
\ba\label{bound_M}
\EE\sup_{t\in [0,T]}|M^\delta_{t}|^2&\lqq 4\EE\int_{0}^{T}\int_{\RR} \Big[H^\delta(Y(s-)+(c_1(X_1(s),v)-c_2(X_2(s),v))) -H^\delta(Y(s-))\Big]^{2}\Pi_{0}(dv)ds.
\ea
Using the same separation argument for (\ref{eq: separation}) 
as in the estimate (\ref{estm: pre H}) we obtain 
\begin{equation}\label{eq: marginale estimate}
\EE\big[\sup_{t\in [0,T]}|M^0_{t}|^2\big]\lqq 4\EE\int_{0}^{T}\int_{\RR}\left(\left|c_1(X_1(s),v)-c_2(X_2(s),v)\right|^{2}\wedge 1\right)\Pi_{0}(dv)ds.
\end{equation}
Using the fact that $|z|^{2}\wedge 1\lqq|z|\wedge 1$ for $z\in\RR$ we get analogously to (\ref{pre Gronwall 1})
\ba
\EE[\sup_{t\in [0,T]}|M_{t}^0|^2]\lqq 16 L_{\Pi}\int_{0}^{T}\EE[(|Y(s)|\wedge1)]ds+4 \Delta. \ea
Inserting (\ref{gronwall_first_mom}) in (\ref{eq: marginale estimate}) and finally both of them 
in (\ref{pre sup}) we conclude 
\ba
&\EE\big[\sup_{t\in[0, T]}\left(|Y(s)|\wedge 1\right)\big]\\
&\lqq  \sqrt{(16L_{\Pi} Te^{(2L_a + 4 L_{\Pi})T}+4)\big(2 T \|a_1 - a_2\|_\infty + 4T\sup_{x\in\RR}\mathcal T_1(\Pi_{1}(x,\cdot),\Pi_{2}(x,\cdot))+2\rho(x_1, x_2)}\big)\\
&\qquad+((2L_a + 4 L_{\Pi}) T e^{(2L_a + 4 L_{\Pi})T}+1)\big(2 T \|a_1 - a_2\|_\infty + 4T\sup_{x\in\RR}\mathcal T_1(\Pi_{1}(x,\cdot),\Pi_{2}(x,\cdot))+2\rho(x_1, x_2)\big).
\ea
Taking $K$ as the maximum of all appearing constants and $x+\sqrt{x} \lqq 2 G(x)$ 
we have obtained the required estimate. This finishes the proof.
\end{proof}

\begin{proof} \emph{of the Proposition \ref{prop:1}: }
Analogously to (\ref{ex:powerlaw}) we get the functions $c_{i}$
\ba
c_{i}(y)=-\left(|y|\lambda_{-}\right)^{1/\alpha^-_{i}}\mathbf{1}_{\{y\lqq -\frac{1}{\lambda_{-}}\}}
+\left(y\lambda_{+}\right)^{1/\alpha^+_i}\mathbf{1}_{\{y\gqq \frac{1}{\lambda_{+}}\}}.
\ea
Thus, the distance $\mathcal T_1$ is calculated as
\ba\label{T_for_polin_for_alpha}
\mathcal T_1(\Pi_{1},\Pi_{2})&=\int_{0}^{\infty}
\left(\mathbf{1}_{\{y\gqq \frac{1}{\lambda_{+}}\}}\left|(y\lambda_{+})^{1/\alpha^+_{1}}-(y\lambda_{+})^{1/\alpha^+_{2}}\right|\wedge 1\right)\frac{dy}{y^{2}}\\
&\qquad +\int_{-\infty}^{0}
\left(\mathbf{1}_{\{y\lqq -\frac{1}{\lambda_{-}}\}}\left||y\lambda_{-}|^{1/\alpha^-_{1}}-|y\lambda_{-}|^{1/\alpha^-_{2}}\right|\wedge 1\right)\frac{dy}{y^{2}}.
\ea
The integrands are integrable with the absolute value since $\alpha_i^\pm>2$.
The first summand equals
\ba
\int_{\frac{1}{\lambda_{+}}}^{\infty}
\left(\left|(y\lambda_{+})^{1/\alpha_{1}}-(y\lambda_{+})^{1/\alpha_{2}}\right|\wedge 1\right)\frac{dy}{y^{2}}
=\la_+\int_{0}^{1}\left(t^{-1/\alpha_{1}}-t^{-1/\alpha_{2}}\right)dt\\
\lqq \frac{\alpha_{1}}{\alpha_{1}-1}-\frac{\alpha_{2}}{\alpha_{2}-1}\lqq \alpha_{2}-\alpha_{1},
\ea
with the help of the change of variables $t=1/(y\lambda_{+})$ and assuming that $2< \alpha_{1}<\alpha_{2}$,
where we have dropped the $+$ superscript for convenience.
The analogous computation for the second integral
in (\ref{T_for_polin_for_alpha}) yields the estimate in Proposition (\ref{prop:1}).
\end{proof}

\end{document}